\documentstyle[12pt]{article}
\roman{enumi} \topmargin=-0.5cm
\title{EXPECTED  NUMBER OF SLOPE CROSSINGS OF  CERTAIN GAUSSIAN   RANDOM  POLYNOMIALS}
\author{S. Rezakhah,  S. Shemehsavar\thanks{ \scriptsize Faculty of Mathematics and Computer Science, Amirkabir University of Technology, 424
 Hafez Avenue, Tehran Iran. email: Rezakhah@aut.ac.ir,$\;\;$ Shemehsavar@aut.ac.ir.}}
\date{}

\begin{document}

\maketitle

\pagestyle{plain}

\renewcommand{\theequation}{\arabic{section}.\arabic{equation}}
\begin{abstract}
 Let $Q_n(x)=\sum _{i=0}^{n} A_{i}x^{i}$ be a random polynomial where the coefficients
 $A_0,A_1,\cdots $ form  a sequence of centered Gaussian random variables. Moreover,
  assume that
the increments $\Delta_j=A_j-A_{j-1}$, $j=0,1,2,\cdots$ are
independent, assuming   $A_{-1}=0$. The coefficients can be
considered as $n$ consecutive observations of a Brownian motion.
We  study the
 number of
 times that such a
random polynomial crosses a line which is not necessarily parallel
to the $x$-axis. More precisely we obtain  the asymptotic behavior
of the expected number of real roots of the equation $Q_n(x)=Kx$,
for the cases that $K$ is any non-zero real constant $K=o(
n^{1/4})$, and $K=o( n^{1/2})$ separately.

\end{abstract}

AMS(2000) subject classifications. Primary 60H42, Secondary 60G99.

 \setcounter{equation}{0}
\newcommand{\be}{\begin{equation}} \newcommand{\ee}{\end{equation}}
\section{Preliminaries}
 \hspace{.2in}

The theory of the expected number of real zeros of random
algebraic polynomials was addressed in  the fundamental work of M.
Kac\cite{kac1} (1943).
 The works of Logan and Shepp \cite{lsh1,lsh2}, Ibragimov and
 Maslova \cite{ibm}, Wilkins \cite{wil},  Farahmand
 \cite{far1}, and Sambandham \cite{sam1,sam2} are other fundamental
 contributions  to the subject. For various aspects on random
 polynomials see Bharucha-Reid and Sambandham \cite{bhs}, and
 Farahmand\cite{far2}.

 There has been recent interest in cases where the
coefficients form certain random processes,
  Rezakhah and Soltani \cite{rs1,rs2}, Rezakhah and Shemehsavar
  \cite{rs3}.

Let  $A_{0},A_{1},\cdots  $ be a mean zero Gaussian random
sequence for which the increments  $\Delta _i =A_i -A_{i-1}, \;
i=1,2,\cdots $ are independent, assuming $A_{-1}=0$. The sequence
$A_0,A_1 \cdots $
 may be considered as successive Brownian points,
 i.e., $A_j=W(t_j),\; j=0,1,\cdots $, where $t_0<t_1< \cdots $
and $\{W(t),\;t\geq 0\}$ is the standard Brownian motion.
In this physical interpretation, Var$(\Delta_j)$ is the distance
 between successive times $t_{j-1},\; t_j$.
  Let
$$ Q_{n}(x)=\sum_{i=0}^{n}A_{i}x^{i},\, \; -\infty<x<\infty . \eqno{(1.1)}$$
We study
 the number of times that $Q_n(x)$ crosses a line which is not necessarily parallel
 to the $x$-axis. Let  $N_{K}(a,b)$
 denote  the  number of real roots of the equation $Q_n(x)=Kx$,
   in the interval $(a,b)$, where
 $K$ is any real constant independent of $x$,
  and multiple roots are counted only once.

 We note that $A_j=\Delta_0 +\Delta _1
+\cdots +\Delta_j, \;\;\; j=0,1,\cdots $, where $ \Delta_i \sim
N(0,\sigma_i^2 )$ and $\Delta_i$ are independent, $i=0,1, \cdots
$.
  Thus $
Q_n(x)=\sum_{k=0}^n (\sum _{j=k}^n x^j)\Delta_k=\sum_{k=0}^n a_k(x)\Delta _k,
$
and $
  Q'_n(x)=\sum_{k=0}^n (\sum _{j=k}^n jx^{j-1})\Delta_k=\sum_{k=0}^n b_k (x)\Delta _k,
$
where
    $$a_k (x) = \sum _{j=k}^n x^j,\;\;
 \; \; \; b_k (x) = \sum _{j=k}^n
jx^{j-1},\;\;\hspace{.3in}  k=0,\cdots , n.\eqno{(1.2)}$$

  Using the function
 $\; \mbox{erf} (t)=2\Phi(t\sqrt{2})-1,\;$ and the fact that
 $t(\mbox{erf}(t))=-t(\mbox{erf}(-t))$,     the classical result of
 Cramer and Leadbetter
[1967 p 285]   can be applied to obtain ;
$$EN_K(a,b)=\int_a^b f_n(x)dx\eqno{(1.3)}$$
where $$
f_n(x)=g_{1}(x)\exp(g_{2}(x))+g_{3}(x)\exp(g_{4}(x))\,\mbox{erf}(g_{5}(x)),
\eqno{(1.4)}$$
$$\hspace{-.2in}\mbox{in which}\hspace{.8in}\;\;\;g_{1}(x)=
\frac{1}{\pi} E A^{-2},\;\;
g_{2}(x)=-\frac{K^2(A^2-2Cx+B^2x^2)}{2E^2},\;\;$$
$$
 g_{3}(x)=\frac{1}{\sqrt{2\pi}} K(A^2-Cx)A^{-3}\;\; g_{4}(x)=
\frac{-K^2x^2}{2A^2},\;\; g_{5}(x)= \frac{K(A^2-Cx)}{AE\sqrt{2} },
 $$
and
$$\;
A^2=\mbox{Var}(Q_n(x))\;=\;\sum_{k=0}^na_k^2(x)\sigma^2_k,\;\;\;\;
B^2=\mbox{Var}(Q'_n(x))=\sum_{k=0}^nb_k^2(x)\sigma^2_k,\;$$
$$
C\!=\!\mbox{Cov}(Q_n(x), Q'_n(x))\!=\!\sum_{k=0}^n
a_k(x)b_k(x)\sigma^2_k,  \;\;\;\; \mbox{and} \;\;\;\;
E^2=A^2B^2-C^2,$$ where  $a_k(x),\;\;$$\;b_k(x)$ is defined by
(1.2).

\section{Asymptotic behaviour of $EN_{K}$ }
 \noindent
 In this section we obtain  the asymptotic behaviour of the expected number of
 real zeros of the equation $Q_{n}(x)=Kx$,  where $Q_{n}(x)$ is
 defined by (1.1).  We prove  the following theorem for the case that
 the increments
 $\Delta_1 \cdots \Delta_n$ have the same distributions, and
 $\sigma_k^2=1$, $k=1\cdots n$.\\ \vspace{-.1in}

 \noindent
 {\bf Theorem 2.1.}   Let $Q_n(x)$ be the random
algebraic polynomial given by (1.1) for which $A_j=\Delta_1+\cdots
\Delta_j$, where $\Delta_i, \;\; i=1,\cdots n$, are independent
and  $\Delta_j\sim N(0,1)$, for $j=1,2,\cdots ,n$. The expected
number of real roots of  $Q_n(x)=Kx$ satisfies the
 following equations:\vspace{.1in}

\noindent (i) -$\;\;\;$
   for  $K=o(n^{1/4})$, and for large $n$\vspace{-.25in}

\begin{eqnarray*}EN_K(-\infty,\infty)\!&\!=\!&\!\frac{1}{\pi}(\log
(2n+1))+\frac{1}{\pi} (1.920134478) \\&&\!\!\!\!-
  \frac{1}{\pi \sqrt{2n}} \left(  \pi -2
\arctan ( \frac {1}{2\sqrt{2n}})\right)
+\!\frac{K^2}{n\pi} (3.126508929)                                   \\
&&\!\!\!+\frac{1}{n\pi}\,C_{1} +o(n^{-1})\hspace{2.5in}(2.1)
\end{eqnarray*}
where  for $n$ odd $C_{1}=1.715215531 ,\;$  and
 for $n$ even  $C_{1}=-0.7200279388.$\vspace{.1in}

\noindent (i$\!\,$i) - $\;\;$  for  $K=o(n^{1/2})$, and for large
$n$\vspace{-.25in}

\begin{eqnarray*}EN_K(-\infty,\infty)=\frac{1}{\pi}(\log
(2n+1))+\frac{1}{\pi} (1.920134478)+o(1)\hspace{1in}(2.2)
\end{eqnarray*}

 \noindent
 {\bf Proof.}
 Due to the behaviour of $Q_n(x)$, the asymptotic behaviour is treated
separately on
  the intervals $1<x<\infty$,
$-\infty<x<-1$, $0<x<1$ and $-1<x<0$.    For $1<x<\infty$, we use
  the change of variable $x=1+\frac{t}{n}$ and
 the equality  $\left( 1+\frac{t}{n}
\right)^n=e^t\left(1-\frac{t^2}{n} \right) +O\left(\frac{1}{n^2}
\right).  $  Using (1.3), we find  that $$EN_K(1 ,\infty
)=\frac{1}{n} \int_{0}^{\infty}f_n\big(1+\frac{t}{n}\big)dt, $$
where  by (1.4) and by tedious manipulation we have that
$$g_{2}\left(1+\frac{t}{n}\right)=o(n^{-2}),\;\;\;g_{3}\left(1+\frac{t}{n}\right)
=o(1),$$ $$ g_{4}\left(1+\frac{t}{n}\right)=o(n^{-2}),
\;\;\;g_{5}\left(1+\frac{t}{n}\right)=o(n^{-1}) \eqno{(2.3)}$$
 and
$$
n^{-1}g_{1}\left(1+\frac{t}{n}\right)=
\frac{1}{\pi}\left(R_1(t)+\frac{S_1(t)}{n}+O\left(
\frac{1}{n^2}\right)\right),\;\;\; n\rightarrow \infty ,
\eqno(2.4)
$$
where
{\small
$$
R_1(t)=\frac
{\sqrt{(4t\!-\!15)e^{4\,t}\!+\!(24t+32)e^{3\,t}\!-\!(8t^
3\!+\!12t^2\!+\!36t\!+\!18)e^{2\,t} +8e^{t}t+1}}{2\,t\left
(-1-3\,{e^{2\,t}}+4\,{e^ {t}}+2\,t{e^{2\,t}}\right )}
$$
and $S_1(t)=S_{11}(t)/S_{12}(t)$ in which
$$
S_{11}(t)=- 0.25\left((4t^2-6t-27)e^{6t}+(156-84t+116t^2-24t^3)e^{5t}
\right.\hspace{1.4in}$$
$$
+( 16t^5- 72t^4+ 96t^3- 212t^2+ 220t- 331)e^{ 4t}
+(328-168t+128t^2-104t^3)e^{3t}
$$
$$\hspace{1in}\left. +( 8t^4 + 8t^3-32t^2+ 42t- 153 )e^{2t}+(28-4t-4t^2)e^t -1\right)
$$
$$
S_{12}(t)=\left((2t-3)e^{ 2t}+ 4e^{t}-
 1\right)^{2}\left((4t-15)e^{4t}+(32+24t){e^{ 3t}}\right.
$$
$$
\hspace{1in}\left. -(8t^3+12t^2+ 36t+ 18)e^{2t} +8te^t  +
1\right)^{1/2}.
$$}
One can easily verify that as $t\rightarrow  \infty$,
$$
R_1(t)=  \frac{1}{2t^{3/2}}+O(t^{-2}),\;\;\;\; \hspace{1in}
S_1(t)= - \frac{1}{8t^{1/2}}+O(t^{-3/2}).
$$
As  (2.4) can not be integrated term by term,
  by noting that
$$\hspace{-1.5in}
\frac{I_{[t>1]}}{8n\sqrt{t}}=\frac{I_{[t>1]}}{8n\sqrt{t}+t\sqrt{t}}
+O\left(\frac{1}{n^2}\right),\hspace{-1in}\eqno(2.5)\hspace{.3in}
$$where
$$I_{[t>1]}=\left\{\begin{array}{ll}1\;\;\;\; \mbox{if }\;\;t\geq 1\\0\;\;\;\;\mbox{if} \;\; t<1\end{array}\right. ,
$$
Also by (2.3), $\; \mbox{erf}(\, g_{5}(t)\,)=o(n^{-1})$. Therefore
 $$
 n^{-1}f_n\bigg(1+\frac{t}{n}\bigg)=
\frac{1}{\pi}\left(R_1(t)+\frac{S_1(t)}{n}\right)+O\left(\frac{1}{n^{2}}\right).
$$
 Thus by (2.5) we have that
$$
n^{-1}f_n\bigg(1+\frac{t}{n}\bigg)+\frac{I_{[t>1]}}{\pi(8n\sqrt{t}+t\sqrt{t})}=\frac{R_1(t)}{\pi}
+\frac{1}{\pi} \left( \frac{S_1(t)}{n}+
\frac{I_{[t>1]}}{8n\sqrt{t}}\right)+O\left(\frac{1}{n^{2}}\right).
$$
This expression is term by term integrable, and provides that
\begin{eqnarray*}EN_K(1,\infty)\!&\!=\!&\! \frac{1}{n}\int_0^{\infty} f_n\bigg(1+\frac{t}{n}\bigg)dt
= \frac{1}{2\pi \sqrt{2n}} \left( - \pi +2
\arctan ( \frac {1}{2\sqrt{2n}})\right)\\
&&\;\;\; +\frac{1}{\pi}\int_0^{\infty}R_1(t)dt+\frac{1}{\pi n}
\int_0^{\infty} \left( S_1(t)+
\frac{I_{[t>1]}}{8\sqrt{t}}\right)dt+O\left(\frac{1}{n^{2}}\right),
\end{eqnarray*}
 where $\;\; \int_0^{\infty}R_1(t)dt\simeq .734874192\;$,
  and $\;\;\int_0^{\infty} \left( S_1(t)+
  \frac{I_{[t>1]}}{8\sqrt{t}}\right)dt=-.25460172372$.

 For  $-\infty<x<-1$, we use  the change of variable
$x=-1-\frac{t}{n}$. So by (1.3) we have that  $EN_K(-\infty, -1
)=\frac{1}{n} \int_{0}^{\infty}
f_n\left(-1-\frac{t}{n}\right)dt.\,$  Using (1.4) we find
$$
n^{-1}g_{1}\left(-1-\frac{t}{n}
\right)=\frac{1}{\pi}\left(R_2(t)+\frac{S_2(t)}{n}
+O\left(\frac{1}{n^2}\right)\right)\eqno(2.6)
$$
where {\small
$$R_2(t)=   1/2\sqrt {{\frac {(1+4t){e^{4\,t}}-(2+4t+12t^2+8t^3){e^{2\,t}}+1}{{t}^{2}
\left( (2t+1){e^
{2\,t}}-1 \right) ^{2}}}}.
$$
 Also we find that for $\,n\,$ even
$\;S_2(t)=(S_{21}(t)+S_{22}(t))/(4S_{23}(t))\;$,
 and for $\,n\,$ odd $\;S_2(t)=
 (S_{21}(t)-S_{22}(t))/(4S_{23}(t))$,
in which
\begin{eqnarray*}
S_{21}(t)\!\!&\!=\!\!&\!  1+\left(
-8\,{t}^{4}+30\,t-8\,{t}^{3}+48\,{t}^{2}-3
 \right) {e^{2\,t}} \\&& + \left( 3-12t+52t^2+96t^3+40t^{4}-16\,{t}^{
5} \right) {e^{4\,t}}\! - \left( 18t+4t^{2}+1 \right) {e^{6\,t}},
\\
S_{22}(t)\!\!&\!=\!\!&\! 4\,{e^{t}}t+ \left(
8\,t+32\,{t}^{3}+40\,{t}^{2}
 \right) {e^{3\,t}}  +
 \left( -8\,{t}^{2}-12\,t \right) {e^{5\,t}}
 \\ S_{23}(t)\!\!&\!=\!\!&\!  \left(e^{4\,t}(4t+1)-2e^{2t}(1+2t+6t^{2}
 +4t^{3})+1\right)^{1/2}
  \left( e^{2\,t}(2t+1)-1
 \right) ^{2}.
\end{eqnarray*}
Also
$$
g_{2}\left(-1-\frac{t}{n}\right)=\left(\frac{K^2}{n}\right)g_{2,1}(t)+
o(n^{-1})\eqno(2.7)
$$
where
$$
g_{2,1}(t)={\frac { \left(
-32\,{t}^{4}-16\,{t}^{3}+16\,{t}^{2}-8\,t \right) {e^{
2\,t}}+8\,t}{ (4t+1){e^{4\,t}}-\left(8\,{t}^{3}
+12\,{t}^{2}+4\,t+2 \right) {e^{2\,t}} +1}}.
$$
 Also we find that
$$
n^{-1}g_{3}\left(-1-\frac{t}{n}\right)=\left(\frac{K}{\sqrt{n\pi}}\right)g_{3,1}(t)+
o(n^{-1})\eqno(2.8)
$$
where
$
g_{3,1}(t)=\frac {- \left(1+\left( 4\,{t}^{2}+2\,t-1 \right)
{e^{2\,t}}\right)}{\sqrt {t} \left( (2t+1)\,{e^{2\,t}}-1 \right)
^{3/2}}.
$
      Also
$$
g_{4}\left(-1-\frac{t}{n}\right)= o\left(n^{-1/2}\right).
\eqno(2.9)
$$
Finally
$$
\hspace{1in}g_{5}\left(-1-\frac{t}{n}\right)=\left(\frac{K}{\sqrt{n}}\right)g_{5,1}(t)+
o(n^{-1}),  \hspace{1.5in} (2.10)
$$
 where \begin{eqnarray*}
 g_{5,1}(t)=- \!\!\! &&\!\!\! 2\sqrt{t} \left(
1+(4t^2+2t-1){e^{2\,t}}
 \right)  \left\{ \left( 6\,t+1+8\,{t}^{2} \right) {e^{6\,t} }\right.
 \\&&+\left( -12\,t-3-32\,{t}^{3}-20\,{
t}^{2}-16\,{t}^{4} \right) {e^{4\,t}}\\&& + \left(
3+6\,t+12\,{t}^{2}+8\,{t }^{3} \right) {e^{2\,t}} \left. -1
\right\}^{-1/2}.
\end{eqnarray*}
  It  can be seen that as $n \rightarrow \infty$,
$$\hspace{-.3in}
R_2(t)=  \frac{1}{2t^{3/2}}+O(t^{-2}),\;\;\;\; \hspace{.1in}
S_2(t)= \frac{-1}{8t^{1/2}}+O(t^{-3/2}),\;\;\;\;\hspace{.1in}
g_{2,1}(t)= o(e^{-t}) ,$$ $$ g_{3,1}(t)=
O(e^{-t}),\;\;\;\hspace{.6in}  g_{5,1}(t)= O(t^{3/2}e^{-t}).
$$
Now by using (1.4) we find  that $\;
f_{n}(-1-t/n)/n:=I_{21}+I_{22}æ\;
 $
 where
$$
I_{21}=\frac{1}{\pi}\left\{R_{2}(t)+\frac{S_{2}(t)}{n}
+\frac{K^2}{n}R_{2}(t)g_{2,1}(t)\right\} +o(n^{-1})
$$
 and
$$
I_{22}=\frac{1}{\pi}\left\{\frac{2K^2}{n}g_{3,1}(t)g_{5,1}(t)\right\}
+o(n^{-1})
 $$
then by using (2.5) we have  that
\begin{eqnarray*}EN_K(-\infty,-1)\!&\!=\!&\! \frac{1}{n}\int_0^{\infty}
f_n\bigg(-1-\frac{t}{n}\bigg)dt \\&\!\!\! = \!\!&\! \frac{1}{2\pi
\sqrt{2n}} \left( - \pi +2 \arctan ( \frac
{1}{2\sqrt{2n}})\right)+\frac{1}{\pi}\int_0^{\infty}\!R_2(t)dt\\
&&+\frac{1}{\pi n} \int_0^{\infty}\! \bigg( S_2(t)+
\frac{I_{[t>1]}}{8\sqrt{t}}\bigg)dt\\&& +
\frac{K^2}{n{\pi}}\int_0^{\infty}\bigg(R_2(t)g_{2,1}(t)
+2g_{3,1}(t) g_{5,1}(t)\bigg) dt  +o(n^{-1}),
\end{eqnarray*}
where
 $\int_0^{\infty}\! \left( S_2(t)\!+\!
 \frac{I_{[t>1]}}{8\sqrt{t}}\right)\!dt\!
 =-0.0322863,\,$  for $n$ odd,  and
$ \int_0^{\infty}\!
\left(S_2(t)\!+\!\frac{I_{[t>1]}}{8\sqrt{t}}\right)\!dt\!=-0.4677136958\,$
for $n$ even.  Also
  $\;\; \int_0^{\infty}\left(R_2(t)\right)dt=1.09564006,\;\;$
and
$$
\int_0^{\infty}\bigg(R_2(t)g_{2,1}(t)
+2g_{3,1}(t)g_{5,1}(t)\bigg) dt=1.593359902.$$

 For $0<x<1$, let $x=1-\frac{t}{n+t}$, then $EN_K(0,1
)=\left(\frac{n}{(n+t)^2} \right) \int_{0}^{\infty} f_n
\left(1-\frac{t}{n+t}\right)dt,$
 in which
$$
 g_{2}\left(1-\frac{t}{n+t}\right)=o(n^{-2}),\hspace{.2in}
\left(\frac{n}{(n+t)^2}
\right)g_{3}\left(1-\frac{t}{n+t}\right)=o(n^{-1})
$$
$$
g_{4}\left(1-\frac{t}{n+t}\right)=o(n^{-2}),\hspace{.2in}
g_{5}\left(1-\frac{t}{n+t}\right)=o(n^{-1})\eqno(2.11)
$$
and
\begin{eqnarray*}
\left(\frac{n}{(n+t)^2} \right)g_{1}\left(1-\frac{t}{n+t}\right)
\!\!\!\!&\!\!\!\!=\!\!\!\!&\!\!\!\! \left(1-\frac{2t}{n} +O\left(
\frac{1}{n^2} \right)\right) \frac{1}{\pi}\left(R_3(t)
+\frac{S_3(t)}{n}+O\left(\frac{1}{n^2}\right)\right)\\
& =\!&\!\!\!\!\!
\frac{1}{\pi}\left(R_3(t)+\frac{S_3(t)-2tR_3(t)}{n}\right)
+O\left(\frac{1}{n^2}\right),\hspace{.5in}(2.12)
\end{eqnarray*}
where we observe that {\small $R_3(t) \equiv  R_1(-t)$ and
 $S_3(t)=S_{31}(t)/S_{32}(t)$,  in which
 \begin{eqnarray*}
S_{31}(t)\!&\!=\!&\!\left(  \left( -7{t}^{2}-{\frac {69}{2}}t-{
\frac {61}{4}} \right) {e^{-6t}}+ \left( 6{t}^{3}+35t-55{t}^{2}+39
 \right) {e^{-5t}}\right.+\\&&\!\!\!\!\! \!\! \left( 49t\!-\!4{t}^{5}+22{t}^
{4}+91{t}^{2}\!-\!{\frac {63}{4}}\!-\!12{t}^{3} \right)
{e^{-4t}}\!-\! \left( 6{t}^{3}+30+44{t}^{2}+66t \right)
{e^{-3t}}\\&&+\left.\left( {\frac
{35}{2}}t+2{t}^{4}-6{t}^{3}+16{t}^{2}+{\frac {123}{4}} \right)
{e^{-2t}}
 + \left( -
9-t-{t}^{2} \right) {e^{-t}}+3/4  \right),
\end{eqnarray*}
and $ S_{32}(t)\equiv S_{12}(-t)$.  We see that, as $t\rightarrow
\infty$, $\; R_3(t)= (2t)^{-1}+O(t^{-1/2}e^{-t/2}),\;$ and $ \;
S_3(t)= \frac{3}{4}+O(t^2e^{-t}) $.
 Now using the  equality
 $$
\frac{I_{[t>1]}}{2t}-\frac{I_{[t>1]}}{4n+2t}=\frac{I_{[t>1]}}{2t}-
\frac{I_{[t>1]}}{4n}+O\left(\frac{1}{n^2}\right),\; \eqno(2.13)$$
we have
\begin{eqnarray*}
\frac{n}{(n+t)^{2}}f_n\bigg(1-\frac{t}{n+t}\bigg)\!&\!
=\!&\!\frac{1}{\pi}\left(R_3(t)-\frac{I_{[t>1]}}{2t}\right)
+\frac{1}{\pi} \left(
 \frac{I_{[t>1]}}{2t}-\frac{I_{[t>1]}}{4n+2t}\right)
\\&&+\frac{1}{n\pi}\left( S_3(t)-2tR_3(t)+\frac{I_{[t>1]}}{4}\right) +O\left(\frac{1}{n^2}\right),
\end{eqnarray*}
Thus
\begin{eqnarray*}
\frac{n}{(n+t)^{2}}\int_0^{\infty}f_n\bigg(1-\frac{t}{n+t}\bigg)dt\!
&\! =\!&\!
\frac{1}{\pi}\int_0^{\infty}\left(R_3(t)-\frac{I_{[t>1]}}{2t}\right)dt+\frac{
 \log (2n+1)}{2\pi}\hspace{.2in}(2.14)
\\&&+\frac{1}{n\pi}\int_0^{\infty}\left( S_3(t)-2tR_3(t)+\frac{I_{[t>1]}}{4}\right)dt  +O\left(\frac{1}{n^2}\right).
\end{eqnarray*} where
$\int_0^{\infty}(R_3(t)-I_{[t>1]}/2t )dt=-0.28977126$ and
$\int_0^{\infty}( S_3(t)-2tR_3(t)+I_{[t>1]}/4)dt=0.497593957.$

For $-1<x<0$, using the change of variable let $x=-1+t(n+t)^{-1}$,
 we have that  $EN_K(-1,0 )\!=\! n(n+t)^{-2} \int_{0}^{\infty}
f_n\left(-1+\frac{t}{n+t}\right)dt, $ we have
$$
\left( \frac{n}{(n+t)^2}\right )g_{1}(t)=\left(
\frac{n^2}{(n+t)^2}\right)
\frac{1}{\pi}\left(R_4(t)+\frac{S_4(t)}{n}+O\left(\frac{1}{n^2}\right)
\right),\eqno(2.15)
$$
in which {\small $R_4(t)\equiv R_2(-t)$.   For $n$ even
$\;S_4(t)=(S_{41}(t)+S_{42}(t))/(4S_{43}(t))$,
 and for $n$ odd $\;S_4(t)=
 (S_{41}(t)-S_{42}(t))/(4S_{43}(t))$, where
 \begin{eqnarray*}
S_{41}(t)\!\! &\!=\!&\!\! 8\left\{\left( {\frac
{15}{4}}t\!-\!7/2{t}^{2}\!-\!3/8 \right) {e^{-6t}}\right.+\left(
15{t}^{4}\!-\!3/2t\!-\!22{t}^{3}+{ \frac
{9}{8}}+19/2{t}^{2}\!-\!2{t}^{5} \right) {e^{-4t}}\\&&+
\left.\left(-9/4t+6{t}^{2}-3{t}^{3}-{\frac {9}{8}}+{t}^{4}
\right) {e^{-2t}} +3/8\right\}
\end{eqnarray*}
 $S_{42}(t)\equiv S_{22}(-t)$,   and
$ S_{43}(t)\equiv S_{23}(-t).\;$
 Also we have that
 $$g_{2}\left(-1+\frac{t}{n+t} \right)=
 \left(\frac{K^2}{n}\right) g_{_{2,1}}(t)+o(n^{-1})
\eqno(2.16)
$$
where
$$
g_{2,1}(t)={\frac {8t \left( 1+(4t^3-2t^2-2t-1){e^{-2\,t}} \right)
}{(4t-1){e^{-4t}}+(2-4t+12t^2-8t^3){e^{-2t}}-1 }}\cdot
$$
Also
$$ \left( \frac{n}{(n+t)^2}\right
)g_{3}(t)=\left(\frac{K}{\sqrt{n\pi}}\right)g_{_{3,1}}(t)
+o(n^{-1})\eqno(2.17)
$$
where
$
g_{_{3,1}}(t)={\frac {-\left( \left( 4\,{t}^{2}-2\,t+1 \right)
{e^{-2\,t}}-1\right)}{\sqrt {t} \left((2t-1){e^{-2\,t}}+1 \right)
^{3/2}}}.
$
  Also
we have that
$$ g_{_{4}}(t)=o (n^{-1/2}
),\eqno(2.18)
$$
Finally we find that
$$
g_{_{5}}(t)=\frac{K}{\sqrt{n}}\;g_{_{5,1}}(t)+o(n^{-1})
\eqno(2.19)$$
 where
\begin{eqnarray*}
g_{_{5,1}}(t)\!\!&\!=\!\!&\!{\frac {  \sqrt {t}\left(\left(
-8{t}^{2}+4t+2 \right){e^{-2t}}-2\right)}{\sqrt
{(2t-1){e^{-2t}}+1}\sqrt
{(1-4t){e^{-4t}}+(8t^3-12t^2+4t-2){e^{-2t}}+1}}}.
\end{eqnarray*}
 As $t\rightarrow \infty $ we have
$$
R_4(t)=\frac{1}{2t}+O\left({t^{1/2}e^{-t}}\right),\;\;\;\;
\hspace{.1in}
S_{4}(t)=\frac{3}{4}+O\left(te^{-t}\right),\;\;\;\;\hspace{.1in}g_{2,1}(t)=-8t+O(t^4e^{-2t})
$$
$$
 g_{3,1}(t)= t^{-1/2}+O(t^{3/2}e^{-2t}),\;\;\;\hspace{.1in}
g_{5,1}(t)= -2\sqrt{t}+O(t^{5/2}e^{-2t}).
$$
then we have by (1.4) that
$$
\frac{n}{(n+t)^{2}}f_n\bigg(-1+\frac{t}{n+t}\bigg):=I_{41}+I_{42}
$$where
\begin{eqnarray*}
I_{41}\!\!&\!=\!&\!\frac{1}{\pi}\frac{n^2}{(n+t)^{2}}\bigg(R_4(t)+\frac{S_{4}(t)}{n}
+\frac{K^2}{n}R_4(t)g_{_{2,1}}(t)\bigg)+o(n^{-1})
\end{eqnarray*}
and
\begin{eqnarray*}
I_{42}\!\!&\!=\!\!&\!\frac{1}{\pi}\frac{n^2}{(n+t)^{2}}\bigg(\frac{2K^2}{n}g_{_{3,1}}(t)g_{_{5,1}}(t)
 \bigg)+o(n^{-1})
\end{eqnarray*}
Now  using (2.13)
 we have that
\begin{eqnarray*}
EN_K(-1,0)\!\!\!&\!=\!&\!\!\! \frac{n}{(n+t)^2}\int_0^{\infty}
f_n\bigg(-1+\frac{t}{n+t}\bigg)\hspace{1.7in}(2.20)\\&=&\!\!\!\!\frac{1}{\pi}\int_0^{\infty}
\left(R_4(t)-\frac{I_{[t>1]}}{2t}\right)dt\\&&\!\!\!
+\frac{1}{\pi}\int_0^{\infty}\! \left(
\frac{I_{[t>1]}}{2t}\!-\!\frac{I_{[t>1]}}{4n+2t} \right)dt
+\frac{1}{n\pi}\int_0^{\infty}\!\! \left(
S_4(t)\!-\!2tR_4(t)\!+\!\frac{I_{[t>1]}}{4}\right)dt\\&&\!+\!
 \frac{K^2}{n\pi}\int_0^{\infty}\!\!\!
\left(R_4(t)g_{2,1}(t)+2g_{3,1}(t)g_{5,1}(t)\right)dt+
o\big(n^{-1}\big).
\end{eqnarray*}
where
  for $n$
odd
 $\;\int_0^{\infty}( S_4(t)-2tR_4(t)+I_{[t>1]}/4)dt
 =1.499908194,\;$ and
 for $n$ even $\; \int_0^{\infty}(S_4(t)-2tR_4(t)
 +I_{[t>1]}/4)dt\!=\!-0.4999082034.$
 Also
$\; \int_0^{\infty}\left(R_{4}(t)\!-\!
I_{[t>1]}(2t)^{-1}\right)dt\!= 0.3793914850,$  and
 $\;
\int_0^{\infty}\left(R_{4}(t)\,g_{2,1}(t)+2\,g_{3,1}(t)g_{5,1}(t)\right)
dt=1.533149028.\; $
 Thus  we arrive at (2.1), and
the first part of the theorem is proved.\\\vspace{-.1in}

 Now for the proof of the second part of the theorem, that is for the case
 $K=o(\sqrt{n})$  we study the asymptotic behavior of $EN_K(a,b)$  for different intervals
 $(-\infty,-1)$, $(-1,0)$, $(0,1)$ and $(1,\infty)$ separately.

 Let   $1<x<\infty$ by using the change of variable
$x=1+\frac{t}{n}$ and (1.4) and by the result (2.4) we find that
\begin{eqnarray*}
\frac{1}{n}\int_{0}^{\infty}f_{n}(1+\frac{t}{n})dt=\frac{1}{\pi}\int_{0}^{\infty}
R_{1}(t)dt+o(n^{-1})
\end{eqnarray*}
where$\int_{0}^{\infty}R_{1}(t)dt=0.734874192$ .\\
For $-\infty<x<-1$ we use  the change of variable
$x=-1-\frac{t}{n} $.  By  the fact  that $K^2/n=o(n)$, we have
$\lim_{\, n\rightarrow \infty} K^2/n=0$. So
$$exp\{K^2g_{_{2,1}}(-1-t/n)/n \}=1+o(1),\;\;\;\;\;
exp\{K^2g_{_{4,1}}(-1-t/n)/n \}=1+o(1).$$ Therefore the relations
(2.7) and (2.9) implies that  $\;exp\{ g_{_2}(-1-t/n)\}=1+o(1)$
and $\;exp\{ g_{_4}(-1-t/n)\}=1+o(1)$. These and  the relations
(1.4), (2.6), (2.8), and (2.10) implies that
\begin{eqnarray*}
\frac{1}{n}\int_{0}^{\infty}f_{n}(-1-\frac{t}{n})dt=\frac{1}{\pi}\int_{0}^{\infty}
R_{2}(t)dt+o(1)
\end{eqnarray*}
where$\int_{0}^{\infty}R_{2}(t)dt=1.095640061$.\\
For $0<x<1\,$ by using the change of variable $x=1-\frac{t}{n+t}$
and relations (1.4), (2.11), (2.12), (2.13) and by a similar
method as in (2.14) we find that
\begin{eqnarray*}
\frac{n}{(n+t)^2}\int_{0}^{\infty}f_{n}(-1-\frac{t}{n+t})dt=
\frac{1}{\pi}\int_{0}^{\infty}(R_{3}(t)-\frac{I[t>1]}{2t})dt
+\frac{1}{2\pi}\log{(2n+1)}
\end{eqnarray*}
where
$\int_{0}^{\infty}(R_{3}(t)-\frac{I[t>1]}{2t})dt=-0.28977126$.\\
For $-1<x<0$ we use the change of variable $x=-1+\frac{t}{n+t}$.
Using (2.16) and (2.18), by the same reasoning  as in the case
$-\infty<x<-1$, we have  that $\;exp\{
g_{_2}(-1+\frac{t}{n+t})\}=1+o(1)$ and $\;exp\{
g_{_4}(-1+\frac{t}{n+t})\}=1+o(1)$. Thus the relations
(2.15),(2.16),(2.17),(2.18),(2.19), and by a similar method as in
(2.20) we find  that
\begin{eqnarray*}
\frac{n}{(n+t)^2}\int_{0}^{\infty}f_{n}(-1+\frac{t}{n+t})dt=
\frac{1}{\pi}\int_{0}^{\infty}(R_{4}(t)-\frac{I[t>1]}{2t})dt
+\frac{1}{2\pi}\log{(2n+1)}
\end{eqnarray*}
where$\int_{0}^{\infty}(R_{4}(t)-\frac{I[t>1]}{2t})dt=0.3793914850$.
 Thus we have $(2.2)$ and the theorem is proved.

\end{document}